\input amstex\documentstyle{amsppt}  
\pagewidth{12.5cm}\pageheight{19cm}\magnification\magstep1
\topmatter
\title Distinguished strata in a reductive group\endtitle
\author G. Lusztig\endauthor
\address{Department of Mathematics, M.I.T., Cambridge, MA 02139}\endaddress
\thanks{Supported by NSF grant DMS-1855773 and by a Simons Fellowship}\endthanks
\endtopmatter   
\document

\define\Irr{\text{\rm Irr}}

\define\si{\sim}

\define\sqc{\sqcup}

\define\qua{\quad}

\define\part{\partial}
\define\emp{\emptyset}

\define\m{\mapsto}
\define\do{\dots}

\define\lra{\leftrightarrow}

\define\sub{\subset}    

\define\T{\times}
\define\ti{\tilde}
\define\nl{\newline}
\redefine\i{^{-1}}

\define\g{\gamma}

\define\e{\epsilon}

\define\io{\iota}

\define\ps{\psi}
\define\r{\rho}
\define\s{\sigma}

\define\k{\kappa}
\redefine\l{\lambda}

\define\Ph{\Phi}

\define\CC{\bold C}

\define\NN{\bold N}

\define\QQ{\bold Q}

\define\SS{\bold S}

\define\ZZ{\bold Z}

\define\ce{\Cal E}

\define\cp{\Cal P}

\define\cs{\Cal S}

\define\cu{\Cal U}

\define\cx{\Cal X}

\define\tA{\ti A}

\define\sha{\sharp}

\define\br{\bar r}

\define\bul{\bullet}

\head Introduction\endhead
\subhead 0.1\endsubhead
Let $G$ be a connected reductive algebraic group over $\CC$, with Weyl group $W$.
In \cite{L15} we have defined a partition of $G$ into finitely many strata. The set of strata
of $G$ is in natural bijection with the image of a map $cl(W)@>>>\Irr(W)$ where
$cl(W)$ is the set  of conjugacy classes in $W$ and $\Irr(W)$ is a set of representatives
for the isomorphism classes of irreducible representations of $W$ over $\QQ$.
We define a partition of $cl(W)$ into subsets (called the {\it strata of $cl(W)$}): a stratum of $cl(W)$ is by definition a nonempty fibre of the map $cl(W)@>>>\Irr(W)$.
We define a partition of $W$ into subsets (called the {\it strata of $W$}): a stratum of $W$
is by definition the inverse image of a stratum of $cl(W)$ under the surjective
map $W@>>>cl(W)$ which takes an element of $W$ to its conjugacy class. The set of strata of
$W$ is in obvious bijection with the set of strata of $cl(W)$ which is in bijection
with the set of strata of $G$.

In this paper we are interested in two problems:

(i) How to parametrize the strata of $cl(W)$ (or $W$)?

(ii) How to describe explicitly each individual stratum of $cl(W)$.
\nl
These problems are solved in \cite{L15}, but we would like to get a simpler are more direct
approach. We shall reduce these problems to the same problems
restricted to a much smaller set of cases.

The set of strata can be viewed as an enlargement of the set of unipotent classes of $G$
(a unipotent class of $G$ is contained in exactly one stratum); this enlargement is built
from the sets of unipotent classes in groups like $G$ but in all characteristics.
According to \cite{BC76}, the classification of unipotent classes of $G$ can be reduced to
the classification of a smaller
set of unipotent classes, the distinguished ones. We will show that, similarly, the
classification of strata of $G$ can be reduced to the classification of a much
smaller set, that of distinguished strata. It turns out that the distinguished strata of $G$
are indexed by a subset $CL_{dist}(W)$ of the set of elliptic conjugacy classes of $W$
which may be called distinguished conjugacy classes. (It would be interesting to find
a description of distinguished conjugacy classes which is purely in terms of $W$ and is
not case by case.) One of our results is that

(a) the set of strata of $G$ (or $W$) has a simple parametrization in terms of the sets
$CL_{dist}(W)$ where $W$ is replaced by various parabolic subgroups of $W$ (see 0.6(a))
to which ``parabolic inclusion'' (see 0.6) is applied.
\nl
Another result of this paper is that

(b) the strata of $cl(W)$ are precisely the connected components of an (oriented) graph with
set of vertices $cl(W)$ (the graph structure on $cl(W)$ is defined in \S2); the
connected components of this graph are remarkably simple: they are products of
Coxeter graphs of type $A$.
\nl
Our definition of the edges of the graph is by first defining (case by case) the edges
for which one end is an elliptic conjugacy class and then applying ``parabolic inclusion''
to such elementary edges where $W$ is replaced by various parabolic subgroups of $W$.
It would be interesting to find a description of these elementary edges which is not case by
case.

In \S3 we show how

(c) $CL_{dist}(W)$ can be parametrized in terms of certain reflection subgroups of $W$.
\nl
Note that by combining (a),(b),(c) one can hope to understand the classification of conjugacy
classes in $W$ in different terms than those in the classification of Carter \cite{C72}.
Namely, the subset $CL_{dist}(W)$ is classified by (c); next, the set of
strata of $cl(W)$ is obtained from (a) by parabolic inclusion and finally the objects
of $cl(W)$ should be described by their position in the graph (product of graphs of type $A$
in (b)) associated to a stratum of $cl(W)$.

We note that the set $CL_{dist}(W)$ and (a),(b) above depend only on $W$ as a Coxeter group
(not on $G$); however, this is not so for (c).

\subhead 0.2\endsubhead
Let $\cp=\{2,3,5,\do\}$ be the set of prime numbers. For any $r\in\cp$ let $G^{(r)}$
be a connected reductive group over an algebraically closed field of characteristic $r$
of the same type as $G$. We set $G^{(0)}=G$. For $r\in\{0\}\cup\cp$ let
$\cu^{(r)}$ be the set of unipotent classes of $G^{(r)}$. By the Springer correspondence
(extended in \cite{L84} to small characteristic) there is a natural imbedding
$\io^{(r)}:\cu^{(r)}@>>>\Irr(W)$ 
whose image is denoted by $\SS^{(r)}(W)$; it is known
that $\SS^{(0)}(W)\sub\SS^{(r)}(W)$.
Let
$$\SS(W)=\cup_{r\in\{0\}\cup\cp}\SS^{(r)}(W)=\cup_{r\in\cp}\SS^{(r)}(W)\sub\Irr(W).$$ 

In \cite{L15} (where the notation $\cs_2(W)$ is used instead of $\SS(W)$)
it is shown that $\SS(W)$ depends only on $W$   
as a Coxeter group, not on the underlying root datum (but it is not clear whether $\SS(W)$ makes sense
for a finite non-crystallographic Coxeter group).

In \cite{L15} we have defined for any $r\in\{0\}\cup\cp$ a surjective map $\k^{(r)}:G^{(r)}@>>>\SS(W)$
whose fibres are called the {\it strata} of $G^{(r)}$; each stratum is a union of conjugacy classes of the
same dimension, independent of $r$ and, according to \cite{C20}, is locally closed in $G^{(r)}$.
If $g\in G^{(r)}$ is unipotent, then $\k^{(r)}(g)$ is the same as the image of the conjugacy class of $g$
under $\io^{(r)}$. It follows that for any $E\in\SS(W)$, there exists $r\in\{0\}\cup\cp$ such that 
the stratum $(\k^{(r)})\i(E)$ contains some unipotent element.

\subhead 0.3\endsubhead
An element of $G^{(r)}$ is said to be distinguished
if it is not contained in a Levi subgroup of a proper parabolic subgroup of $G^{(r)}$ (see \cite{BC76}).
Let $\cu^{(r)}_{dist}$ be the set of unipotent classes in $G^{(r)}$ in which some/any element is
distinguished. Let $\SS^{(r)}_{dist}(W)=\io^{(r)}(\cu^{(r)}_{dist})$.

We say that $E\in\SS(W)$ is distinguished if $E\in\cup_{r\in\{0\}\cup\cp}\SS^{(r)}_{dist}(W)$
or equivalently if there exists $r\in\{0\}\cup\cp$ such that 
the stratum $(\k^{(r)})\i(E)$ contains some distinguished unipotent element.
Let $\SS_{dist}(W)$ be the set of distinguished elements of $\SS(W)$.

In \S1 we will show:

(a) {\it $E\in\SS(W)$ is distinguished if and only if there exists $r\in\{0\}\cup\cp$ such that the
stratum $(\k^{(r)})\i(E)$ contains some distinguished (not necessarily unipotent) element of $G^{(r)}$.} 
\nl
A stratum of $G^{(r)}$ (with $r\in\{0\}\cup\cp$)
is said to be distinguished if it is of the form $(\k^{(r)})\i(E)$ where
$E\in\SS_{dist}(W)$ (such a stratum need not contain a distinguished unipotent element).

\subhead 0.4\endsubhead
For $C\in cl(W)$ let $m(C)$ be the dimension of the $1$-eigenspace of some/any $w\in C$ on the
reflection representation of $W$.
We shall write $\Ph:cl(W)@>>>\SS(W)$ for what in \cite{L15} is denoted by ${}'\Ph$ (a surjective
map).
In \cite{L15} it is shown that

(a) {\it for any $E\in\SS(W)$ there is a unique $C_E\in\Ph\i(E)$ which is as elliptic as
possible, that is $m(C_E)\le m(C)$ for any $C\in\Ph\i(E)$;}
\nl
thus $E\m C_E$ is a cross section of the surjective map $\Ph$. The following variant
of (a) will be verified in \S1.

(b) {\it for any $E\in\SS(W)$ there is a unique $C'_E\in\Ph\i(E)$ which is as non-elliptic as possible,
that is $m(C'_E)\ge m(C)$ for any $C\in\Ph\i(E)$.}
\nl
Let $CL(W)$ be the image of the map $E\m C'_E$, $\SS(W)@>>>cl(W)$. Note that $\Ph$ restricts to a bijection
$CL(W)@>\si>>\SS(W)$. Under this bijection, the subset $\SS_{disc}(W)$ of $\SS(W)$ corresponds to
a subset $CL_{dist}(W)$ of $CL(W)$. The conjugacy classes of $W$ contained in $CL_{dist}(W)$
are said to be
distinguished. The following result will be proved in \S1.

(c) {\it Let $C\in CL(W)$. We have $C\in CL_{dist}(W)$ if and only if $C$ is elliptic (that is, $m(C)=0$).}

\subhead 0.5 \endsubhead
For $r\in\{0\}\cup\cp$ let
$\ps^{(r)}:cl(W)@>>>\cu^{(r)}$ be the surjective map
defined in \cite{L11a}. (In the case $r=0$, an alternative definition of this map is
given in \cite{Y20}.) Let $\Ph^{(r)}:cl(W)@>>>\Irr(W)$ be the composition of $\ps^{(r)}$
with $\io^{(r)}:\cu^{(r)}@>>>\Irr(W)$.
From the explicit description of $\ps^{(r)}$ in \cite{L12} and the explicit description
of $CL_{dist}(W)$ given in this paper we see that

(a) {\it if $C\in CL_{dist}(W)$ then $\Ph^{(r)}(C)$ is independent of $r$. Hence, by the definition
of $\Ph$ in \cite{L15, 4.1}, we have $\Ph(C)=\Ph^{(r)}(C)$ for all $r$.}

\subhead 0.6 \endsubhead
Let $\{s_i;i\in I\}$ be the set of simple reflections of $W$. For $J\sub I$ let $W_J$
be the subgroup of $W$ generated by $\{s_i;i\in J\}$; this is the Weyl group of a Levi
subgroup of a parabolic subgroup of $G$. Hence $CL(W_J)$ and its subset $CL_{dist}(W_J)$
are defined. For $C_1\in cl(W_J)$ we define $\r_J(C_1)\in cl(W)$ by the condition
$C_1\sub\r_J(C_1)$. Now $C_1\m\r_J(C_1)$ is an injective map $cl(W_J)@>>>cl(W)$.
(We call it {\it parabolic inclusion}.)
If $J\sub I,J'\sub I$, we say that $J,J'$ are equivalent if
$W_J,W_{J'}$ are conjugate under an element of $W$. Let $\cx$ be a set of representatives
for the equivalence classes of subsets $J\sub I$ for the equivalence relation above.
The following result can be deduced from the explicit description of $CL(W)$ given in
\cite{L15} and that of $CL_{dist}(W)$ given in this paper.

(a) $CL(W)=\sqc_{J\in\cx}\r_J(CL_{dist}(W_J))$.

\subhead 0.7. Notation \endsubhead
A bipartition is a sequence $\l_*=(\l_1,\l_2,\l_3,\do)$ in $\NN$ such that $\l_i=0$ for
$i$ large and $\l_1\ge\l_3\ge\l_5\ge\do$, $\l_2\ge\l_4\ge\l_6\ge\do$. We write
$|\l_*|=\l_1+\l_2+\l_3+\do$. Let $BP$ be the set of bipartitions.

Let $T$ be the set of all $\l_*\in BP$ such that $\l_1\ge\l_2\ge\l_3\ge\l_4\ge\do$.
Let $R$ be the set of all $\l_*\in T$ such that $\l_i$ is even for any $i$.
Let $P$ be the set of all $\l_*\in T$ such that $\l_1=\l_2,\l_3=\l_4,\l_5=\l_6,\do$,
For $\l_*\in T$ and $j>0$ we set $\mu_j(\l_*)=\sha(\{k\ge1;\l_k=j\})$.
For $m\in\NN$ let $BP^m=\{\l_*\in BP,|\l_*|=m\}$, $T^m=T\cap BP^m$.
Let $T_{ev}$ (resp. $R_{ev}$) be the subset of $T$ (resp. $R$) consisting of the $\l_*$
with an even number of $>0$ terms.
For $a,b$ in $\ZZ$ we write $a\gg b$ instead $a-b\ge2$.

\head 1. Proof of 0.3(a), 0.4(b), 0.4(c)\endhead
\subhead 1.1\endsubhead
In this section we prove 0.3(a), 0.4(b), 0.4(c). To do this we can assume that $G$ is almost simple.
It is also enough to consider only one such $G$ in each isogeny class.
The case where $G$ is of classical (resp. exceptional) type is considered in 1.2-1.9 (resp.
1.10-1.15).

\subhead 1.2\endsubhead
Assume first that $G=SL_n(\CC), n\ge2$. In this case $\SS^{(r)}(W)=\Irr(W)=\SS(W)$ for any $r$
and the map $\Ph$ is a bijection $cl(W)@>\si>>\SS(W)$. In this case 0.4(b) is obvious and we have
$CL(W)=cl(W)$. Also 0.3(a) is immediate (an element
is distinguished if and only if is regular unipotent
times a central element).
Note that $CL_{dist}(W)$ consists of a single element:
the class of the Coxeter element; thus 0.4(c) holds.

\subhead 1.3\endsubhead
Until the end of 1.5 we assume that $G=Sp_{2n}(\CC),n\ge2$. 
Then $\SS(W)=\SS^{(2)}(W)$ and $\Ph$ becomes a map $cl(W)@>>>\SS^{(2)}(W)$.
By \cite{L12} we have bijections

(a) $\cu^{(2)}\lra$ (set of all pairs $(c_*,\e)$ where $c_*\in T^{2n}$ is such that
$\mu_j(c_*)\in2\NN$ for any odd $j$ and
$\e:\{j\in2\NN+2;\mu_j(c_*)\in2\NN+2\}\}@>>>\{0,1\}$),

(b) $cl(W)\lra(R\T P)^{2n}:=\{(r_*,p_*)\in R\T P;|r_*|+|p_*|=2n\}$.

Via (a),(b), $(\io^{(2)})\i\Ph:cl(W)@>>>\cu^{(2)}$ becomes the map

(c) $(r_*,p_*)\m(c_*,\e)$ where $\mu_j(c_*)=\mu_j(r_*)+\mu_j(p_*)$ for $j>0$
and for any $j\in2\NN+2$ such that $\mu_j(c_*)\in2\NN+2$, we have
$\e(j)=1$ if $j=r_i$ for some $i$ and $\e(j)=0$, otherwise.
\nl
Via (b), the map $cl(W)@>>>\NN, C\m m(C)$ becomes the map which to any
$(r_*,p_*)\in(R\T P)^{2n}$ associates $|p_*|/2$.
To prove 0.4(b) we must show that, if $(c_*,\e)$ (as in (a)) is given, then
there is a unique $(r_*,p_*)\in(R\T P)^{2n}$ which maps to it (as in (c))
and has $|p_*|$ maximum possible. Thus,

for $j\in2\NN+2$ such that $\mu_j(c_*)\in2\NN+2$, we must have that
$\mu_j(r_*)=2,\mu_j(p_*)=\mu_j(c_*)-2$ (if $\e(j)=1$) and $\mu_j(r_*)=0,\mu_j(p_*)=\mu_j(c_*)$
(if $\e(j)=0$); for $j\in2\NN+2$ such that $\mu_j(c_*)\in2\NN+1$, we must have that
$\mu_j(r_*)=1,\mu_j(p_*)=\mu_j(c_*)-1$; for $j\in2\NN+2$ such that $\mu_j(c_*)=0$, we must have
that $\mu_j(r_*)=0,\mu_j(p_*)=0$; for $j\in2\NN+1$, we must have that $\mu_j(r_*)=0,\mu_j(p_*)=\mu_j(c_*)$.
\nl
This proves 0.4(b) in our case.

Note that:

(e) {\it $CL(W)$ is (via (b)) the set of all $(r_*,p_*)\in(R\T P)^{2n}$ 
such that $\mu_j(r_*)\le2$ for any $j>0$.}

\subhead 1.4\endsubhead
As in \cite{L15} we have a bijection

(a) $\Irr(W)\lra BP^n$.
\nl
Using \cite{L15, \S3}, we see that 

(b) {\it when $r\ne2$, the subset $\io^{(r)}(\cu^{(r)}_{dist})$ of $\Irr(W)$
becomes via (a) the subset of $BP^n$ consisting of sequences of the form
$(a_1>a_2>\do>a_s>0,0,0,\do)$.}
\nl
By \cite{W63} (see also \cite{LS12, 6.2}),
the set $\cu^{(2)}_{dist}$ can be identified via 1.3(a) with
the subset of $\cu^{(2)}$ consisting of

(c) {\it all $(c_*,\e)$ (as in 1.3(a)) such that $\mu_j(c_*)=0$ for odd $j$, $\mu_j(c_*)\le2$ for
even $j$ and $\e(j)=1$ whenever $j$ is even and $\mu_j(c_*)=2$.}
\nl
Using \cite{L15, \S3} we see that the subset $\io^{(2)}(\cu^{(2)}_{dist})$ of $\Irr(W)$ becomes
via (a) the subset of $BP^n$ formed by the sequences $(c_1/2,c_2/2,c_3/2,\do)$ for various
$(c_*,\e)$ as in (c). This is the same as the set of all $(a_1\ge a_2\ge a_3\ge\do)\in BP^n$
such that there are no consecutive equalities between the non zero $a_i$.
This set contains the set (b). It follows that

(d) $\SS^{(r)}_{dist}(W)\sub\SS^{(2)}_{dist}(W)$ for any $r$.
\nl
Under our bijection $CL_{dist}(W)\lra \SS_{dist}(W)$,
the set of $(c_*,\e)$ as in (c) corresponds to the set of $(r_*,p_*)\in(R\T P)^{2n}$
such that $\mu_j(r_*)\le2$ for $j>0$ and $p_*=(0,0,0,\do)$; this is the same as the set of
all $(r_*,p_*)\in(R\T P)^{2n}$ which under 1.3(b) correspond to
elliptic conjugacy classes in $W$ which are in $CL(W)$. 
This implies (by (d))
that 0.4(c) holds in our case.

\subhead 1.5\endsubhead
Let $g\in G^{(r)}$ be a distinguished element. To prove 0.3(a) it is enough to show
that $\k^{(r)}(g)\in\SS_{dist}(W)$. If $r=2$ then $g$ is unipotent and the result is clear.
Thus we can assume that $r\ne2$. Using \cite{L15, \S3} we see that under the bijection 1.4(a), 
$\k^{(r)}(g)$ corresponds to a bipartition of the form

(a) $((a_1+b_1)/2,(a_2+b_2)/2,(a_3+b_3)/2,\do)$
\nl
where
$$a_1=\nu_1/2,a_2=\nu_2/2,a_3=\nu_3/2,\do,a_s=\nu_s/2,a_{s+1}=0,a_{s+2}=0,\do,$$
$$b_1=\nu'_1/2,b_2=\nu'_2/2,b_3=\nu'_3/2,\do,b_t=\nu'_t/2,b_{t+1}=0,b_{t+2}=0,\do,$$
and $\nu_1>\nu_2>\nu_3>\do>\nu_s$, $\nu'_1>\nu'_2>\nu'_3>\do>\nu'_t$ are even integers $\ge2$
with $\sum_k\nu_k+\sum_k\nu'_k=2n$.

Clearly, $(a_1+b_1)/2\ge(a_2+b_2)/2\ge(a_3+b_3)/2\ge\do$. If 
$(a_i+b_i)/2=(a_{i+2}+b_{i+2})/2$
then $a_i=a_{i+2},b_i=b_{i+2}$ hence $a_{i+2}=0,b_{i+2}=0$ and $(a_{i+2}+b_{i+2})/2=0$.
Thus (a) corresponds to an element of $\SS^{(2)}_{dist}(W)$. This proves 0.3(a) in our case.

\subhead 1.6\endsubhead
We now assume that $G=SO_{2n+1}(\CC), n\ge3$.
Then the arguments in 1.3, 1.4 can be used word by word in the present case except that
1.4(b) must be replaced by the following statement:

(a) When $r\ne2$, the subset $\io^{(r)}(\cu^{(r)}_{dist})$ of $\Irr(W)$
becomes via 1.4(a) the subset of $BP^n$ consisting of sequences $(c_1,c_2,c_3,\do)$ where
$$\align&c_1=(\nu_1-1)/2,c_2=(\nu_2+1)/2,c_3=(\nu_3-1)/2,\\&
c_4=(\nu_4+1)/2,\do,c_{2s+1}=(\nu_{2s+1}-1)/2,c_{2s+2}=0, c_{2s+3}=0,\do\endalign$$
and $\nu_1>\nu_2>\nu_3>\do>\nu_{2s+1}$ are odd integers $\ge1$ with sum $2n+1$.

(Note that $c_1,c_2,\do)\in P$, with no two successive equalities between its nonzero terms
hence it corresponds to an element in $\SS^{(2)}_{dist}(W)$.)

The argument in 1.5 also continues to hold except that 1.5(a) must be replaced by

(b) $((a_1+b_1)/2,(a_2+b_2)/2,(a_3+b_3)/2,\do)$
\nl
where
$$\align&a_1=(\nu_1-1)/2,a_2=(\nu_2+1)/2,a_3=(\nu_3-1)/2,\\&
a_4=(\nu_4+1)/2,\do,a_{2s+1}=(\nu_{2s+1}-1)/2,a_{2s+2}=0, a_{2s+3}=0,\do,\endalign$$
$$\align&b_1=(\nu'_1+1)/2,b_2=(\nu'_2-1)/2,b_3=(\nu'_3+1)/2,\\&
b_4=(\nu'_4-1)/2,\do,b_{2t}=(\nu'_{2t}-1)/2,b_{2t+1}=0, b_{2t+2}=0,\do,\endalign$$
and $\nu_1>\nu_2>\nu_3>\do>\nu_{2s+1}$, $\nu'_1>\nu'_2>\nu'_3>\do>\nu'_{2t}$ are odd integers $\ge1$
with $\sum_k\nu_k+\sum_k\nu'_k=2n+1$.

We have $a_1\ge a_2\ge a_3\ge\do$ and $b_1\ge b_2\ge b_3\ge\do$ hence
$$(a_1+b_1)/2\ge(a_2+b_2)/2\ge(a_3+b_3)/2\ge\do.$$
If $(a_i+b_i)/2=(a_{i+2}+b_{i+2})/2$ then $a_i=a_{i+2},b_i=b_{i+2}$ hence $a_i=b_i=0$.
\nl
Thus (a) corresponds to an element of $\SS^{(2)}_{dist}(W)$. This proves 0.3(a) in our case.

\subhead 1.7\endsubhead
Until the end of 1.9 we assume that $G=SO_{2n}(\CC), n\ge4$.
We have $\SS(W)=\SS^{(2)}(W)$ and $\Ph$ becomes a map $cl(W)@>>>\SS^{(2)}(W)$.
Now each of
$$cl(W),\Irr(W),\SS^{(2)}(W),\cu^{(2)}$$ has a natural involution
induced by conjugation by an element in the non-identity component of the
full orthogonal group. We then have partitions
$$cl(W)=cl(W)'\sqc cl(W)'',\Irr(W)=\Irr(W)'\sqc\Irr(W)'',$$
$$\SS^{(2)}(W)=\SS^{(2)}(W)'\sqc\SS^{(2)}(W)'',\cu^{(2)}=\cu^{(2)}{}'\sqc\cu^{(2)}{}'',$$
where $()'$ denotes the set of fixed points of the involution and $()''$ denotes its
complement.

By \cite{L12} we have bijections

(a)  (set of orbits of the involution on $\cu^{(2)}{}''$) $\lra$
(set of all pairs $(c_*,\e)$ as in 1.3(a) such that $\e=0$ and
$\mu_j(c_*)=0$ for $j$ odd, $\mu_j(c_*)$ even for $j$ even),

(b) (set of orbits of the involution on $cl(W)''$) $\lra$
(set of all pairs $(r_*,p_*)\in(R\T P)^{2n}$
such that $r_*=(0,0,\do)$ and
$\mu_j(p_*)=0$ for $j$ odd),

(c) $\cu^{(2)}{}'\lra$
(set of all pairs $(c_*,\e)$ as in 1.3(a) which are not as in (a) and are such that
$c_*\in T_{ev}$),

(d) $cl(W)'\lra$ (set of all pairs $(r_*,p_*)\in(R\T P)^{2n}$ which are not as in (b)
and are such that $r_*\in R_{ev}$).
\nl
Now $(\io^{(2)})\i\Ph$ induces the (surjective) map $cl(W)'@>>>\cu^{(2)}{}'$ which by (c),(d)
becomes the map $(r_*,p_*)\m(c_*,\e)$ given by the same rule as in 1.3(c).
The same proof as in 1.3 shows that if $(c_*,\e)$ (as in (c)) is given, then
there is a unique $(r_*,p_*)$ (as in (d)) which maps to it and has $|p_*|$ maximum possible.
This implies that 0.4(b) holds for any $E\in\SS(W)=\SS^{(2)}(W)$ that is contained in
$\SS^{(2)}(W)'$. If $E\in\SS^{(2)}(W)''$ then $\Ph\i(E)$ consists of a single element so
that 0.4(b) holds automatically for such $E$. This proves 0.4(b) in our case.

In our case

(e) the set $CL(W)$ becomes the set of pairs $(r_*,p_*)$ as in 1.3(e)
such that $r_*\in R_{ev}$ and with each pair as in (b) repeated twice.

\subhead 1.8\endsubhead
As in 1.4(a) we have a bijection

(a) $\Irr(W)'\lra$ set of all $(a_1,a_2,a_3,\do)\in BP^n$ such that $a_1-a_2,a_3-a_4,\do$ are
not all zero.

Using \cite{L15, \S3} we see that 

(b) when $r\ne2$, the subset $\io^{(r)}(\cu^{(r)}_{dist})$ of $\Irr(W)'$
becomes via (a) the subset of $BP^n$ consisting of sequences of the form
$$((\nu_1+1)/2,(\nu_2-1)/2,(\nu_3+1)/2,\do,(\nu_{2s}-1)/2,0,0,0,\do)$$
where $\nu_1>\nu_2>\do>\nu_{2s}$ are odd $\ge1$.

By \cite{W63} (see also \cite{LS12, 6.2}),
the set $\cu^{(2)}_{dist}$ can be identified via 1.7(a) with
the set consisting of

(c) all $(c_*,\e)$ (as in 1.3(a)) where $c_*\in T_{ev}$ and
such that $\mu_j(c_*)=0$ for odd $j$, $\mu_j(c_*)\le2$ for
even $j$ and $\e(j)=1$ whenever $j$ is even and $\mu_j(c_*)=2$.

Using \cite{L15, \S3} we see that the subset $\io^{(2)}(\cu^{(2)}_{dist})$ of $\Irr(W)'$ becomes
via (a) the subset of $BP^n$ formed by the sequences
$$((c_1+2)/2,(c_2-2)/2,(c_3+2)/2,(c_4-2)/2,\do,(c_{2s}-2)/2,0,0,0,\do)$$
for various $(c_*,\e)$ as in (c) with $c_*=(c_1\ge c_2\ge\do\ge c_{2s}>0,0,0,\do)$.

If $\nu_1>\nu_2>\do>\nu_{2s}$ are odd $\ge1$ then 
$$\align&((\nu_1+1)/2,(\nu_2-1)/2,(\nu_3+1)/2,\do,(\nu_{2s}-1)/2,0,0,0,\do)\\&=
((c_1+2)/2,(c_2-2)/2,(c_3+2)/2,(c_4-2)/2,\do,(c_{2s}-2)/2,0,0,0,\do)\endalign$$
where
$c_1=\nu_1-1, c_2=\nu_2+1,c_3=\nu_3-1,\do, c_{2s}=\nu_{2s}+1$
are even and $c_1\ge c_2\gg c_3\ge c_4\gg\do\ge c_{2s}\gg0$.
From this we see that $\io^{(r)}(\cu^{(r)}_{dist})\sub\io^{(2)}(\cu^{(2)}_{dist})$ for any $r$, that is

(d) $\SS^{(r)}_{dist}(W)\sub\SS^{(2)}_{dist}(W)$ for any $r$.
\nl
Under our bijection $CL_{dist}(W)\lra \SS_{dist}(W)$,
the set of $(c_*,\e)$ as in (c) corresponds to the set of $(r_*,p_*)$
as in 1.7(d) such that $\mu_j(r_*)\le2$ for $j>0$ and $p_*=(0,0,0,\do)$; this is the same as the set of
all $(r_*,p_*)$ in 1.7(e) which correspond to elliptic conjugacy classes in $W$. This implies (by (d))
that 0.4(c) holds in our case.

\subhead 1.9 \endsubhead
Let $g\in G^{(r)}$ be a distinguished element. To prove 0.3(a) it is enough to show
that $\k^{(r)}(g)\in\SS_{dist}(W)$. If $r=2$ then $g$ is unipotent and the result is clear.
Thus we can assume that $r\ne2$. Using \cite{L15, \S3} we see that under the bijection 1.8(a), 
$\k^{(r)}(g)$ corresponds to a bipartition of the form

(a) $((a_1+b_1)/2,(a_2+b_2)/2,(a_3+b_3)/2,\do)$
\nl
where
$$\align &
a_1=(\nu_1+1)/2,a_2=(\nu_2-1)/2,a_3=(\nu_3+1)/2,a_4=(\nu_4-1)/2,\do,\\&a_{2s}=(\nu_{2s}-1)/2,
a_{2s+1}=0, a_{2t+2}=0,\do,\endalign$$
$$\align&b_1=(\nu'_1+1)/2,b_2=(\nu'_2-1)/2,
b_3=(\nu'_3+1)/2,\\&b_4=(\nu'_4-1)/2,\do,b_{2t}=(\nu'_{2t}-1)/2,
b_{2t+1}=0, b_{2t+2}=0,\do,\endalign$$
and
$$\nu_1>\nu_2>\nu_3>\do>\nu_{2s}, \qua \nu'_1>\nu'_2>\nu'_3>\do>\nu'_{2t}$$
are odd integers $\ge1$
with $\sum_k\nu_k+\sum_k\nu'_k=2n$. We can assume that $t\le s$. We have
$$\align&((a_1+b_1)/2,(a_2+b_2)/2,(a_3+b_3)/2,\do)\\&
=((c_1+2)/2,(c_2-2)/2,(c_3+2)/2,(c_4-2)/2,\do,(c_{2s}-2)/2,0,0,0,\do)\endalign$$
where
$$\align&c_1=\nu_1+\nu'_1,c_2=\nu_2+\nu'_2,\do,c_{2t}=\nu_{2t}+\nu'_{2t}, \\&
c_{2t+1}=\nu_{2t+1}-1,c_{2t+2}=\nu_{2t+2}+1,\do, c_{2s}=\nu_{2s}+1.\endalign$$

(The terms $c_{2t+1},\do,c_{2s}$ are missing if $t=s$.)
Note that $c_1,c_2,\do,c_{2s}$ are even, nonzero and
$c_1\ge c_2\ge c_3\ge c_4\ge\do\ge c_{2s}$ with no consecutive equalities.
Thus (a) corresponds to an element of $\SS^{(2)}_{dist}(W)$. This proves 0.3(a) in our case.

\subhead 1.10\endsubhead   
In the remainder of this section we assume that $G$ is simple of exceptional type.
In this case 0.4(b) can be deduced from tables in \cite{L15}.
In 1.11-1.15 we describe in each case, using notation of
Carter \cite{C72} and that of \cite{L15},
the bijection between $CL(W)$ and $\SS(W)$ in the form $C\lra E$.
Here we also use the description of distinguished unipotent elements in \cite{M80}, \cite{LS12}.
Now 0.3(a),0.4(c)  can be verified in each case using the definitions.

\subhead 1.11. Type $G_2$\endsubhead   

$[G_2]\lra1_0$; $[A_2]\lra2_1$; $[A_1+\tA_1]\lra2_2$;

$[\tA_1]\lra1_3$; 

$[A_1]\lra1_3$;  

$[A_0]\lra1_6$.  

Here the first three items are in $CL_{dist}(W)$.
Note that $2_2\in\SS^{(r)}_{dist}(W)$ for a single $r$ namely $r=3$.

\subhead 1.12. Type $F_4$ \endsubhead   

$[F_4]\lra1_0$; $[B_4]\lra4_1$; $[F_4(a_1)]\lra9_2$; $[D_4(a_1)]\lra12_4$;
$[A_3+\tA_1]\lra16_5$; $[\tA_2+\tA_2]\lra6_6$;

$[B_3]\lra8_3$;       

$[C_3]\lra8_3$;         

$[A_3]\lra9_6$;       

$[B_2+A_1]\lra9_6$;

$[A_2+\tA_1]\lra4_7$;

$[\tA_2+A_1]\lra4_7$;

$[B_2]\lra 4_8$;     

$[\tA_2]\lra 8_9$;        

$[A_2]\lra 8_9$; 

$[A_1+\tA_1]\lra9_{10}$;

$[2A_1]\lra 4_{13}$;         

$[A_1]\lra 2_{16}$;       

$[\tA_1]\lra 2_{16}$; 

$[A_0]\lra 1_{24}$.

Here the first six items are in $CL_{dist}(W)$.

Note that $16_5,6_6$ are in $\SS^{(r)}_{dist}(W)$ for a single $r$ namely $r=2$.

\subhead 1.13. Type $E_6$\endsubhead              

$[E_6]\lra1_0$; $[E_6(a_1)]\lra6_1$; $[E_6(a_2)]\lra 30_3$,

$[D_5]\lra20_2$; $[D_5(a_1)]\lra64_4$

$[A_5]\lra15_4$   

$[A_4+A_1]\lra60_5$

$[2A_2+A_1]\lra 10_9$

$[D_4]\lra 24_6$; $[D_4(a_1)]\lra 80_7$

$[A_4]\lra81_6$

$[A_3+A_1]\lra60_8$    

$[A_2+2A_1]\lra 60_{11}$

$[2A_2]\lra 24_{12}$

$[A_3]\lra 81_{10}$

$[A_2+A_1]\lra 64_{13}$

$[3A_1]\lra 15_{16}$

$[A_2]\lra 30_{15}$

$[2A_1]\lra 20_{20}$      

$[A_1]\lra 6_{25}$   

$[A_0]\lra 1_{36}$   

Here the first three items are in $CL_{dist}(W)$.

\subhead 1.14. Type $E_7$\endsubhead 

$[E_7]\lra 1_0$; $[E_7(a_1)]\lra 7_1$; $[E_7(a_2)]\lra 27_2$;
$[E_7(a_3)]\lra 56_3$; $[A_7]\lra189_5$; $[E_7(a_4)]\lra315_7$;

$[E_6]\lra 21_3$; $[E_6(a_1)]\lra120_4$; $[E_6(a_2)]\lra405_8$;

$[D_6]\lra35_4$; $[D_6(a_1)]\lra210_6$; $[D_6(a_2)]\lra280_8$; $[2A_3]\lra378_{14}$;

$[A_6]\lra105_6$;

$[D_5+A_1]\lra168_6$; $[D_5(a_1)+A_1]\lra378_9$;

$[(A_5+A_1)']\lra70_9$;

$[A_4+A_2]\lra210_{10}$;

$[A_3+A_2+A_1]\lra210_{13}$;

$[D_5]\lra189_7$; $[D_5(a_1)]\lra420_{10}$;

$[A''_5]\lra216_9$;

$[A_4+A_1]\lra512_{11}$;

$[A'_5]\lra105_{12}$;

$[D_4+A_1]\lra84_{12}$; $[D_4(a_1)+A_1]\lra405_{15}$;           
                             
$[A_3+A_2]\lra84_{15}$;

$[(A_3+2A_1)']\lra216_{16}$;

$[2A_2+A_1]\lra70_{18}$;

$[A_2+3A_1]\lra105_{21}$;

$[A_4]\lra420_{13}$;

$[D_4]\lra105_{15}$; $[D_4(a_1)]\lra315_{16}$;                  

$[(A_3+A_1)'']\lra280_{17}$;

$[(A_3+A_1)']\lra189_{20}$;

$[2A_2]\lra168_{21}$;

$[A_2+2A_1]\lra189_{22}$;

$[(4A_1)']\lra15_{28}$;

$[A_3]\lra210_{21}$;

$[A_2+A_1]\lra120_{25}$;

$[(3A_1)'']\lra35_{31}$;

$[(3A_1)']\lra21_{36}$;

$[A_2]\lra56_{30}$;

$[2A_1]\lra27_{37}$;

$[A_1]\lra7_{46}$;   

$[A_0]\lra1_{63}$.   

Here the first six items are in $CL_{dist}(W)$.

There are two conjugacy classes $4A_1$; one of them, $(4A_1)''$, comes from $W_J$ of type
$D_4$.
There are two conjugacy classes $A_5+A_1$; one of them, $(A_5+A_1)''$, comes from $W_J$ of type
$E_6$.

\subhead 1.15. Type $E_8$\endsubhead

$[E_8]\lra 1_0$; $[E_8(a_1)]\lra 8_1$; $[E_8(a_2)]\lra 35_2$; $[E_8(a_4)]\lra 112_3$;

$[E_8(a_5)]\lra 210_4$; $[D_8]\lra 560_5$; $[E_8(a_3)]\lra 700_6$; $[E_8(a_7)]\lra1400_7$;

$[E_8(a_6)]\lra1400_8$; $[D_8(a_2)]\lra 3240_9$; $[A_8]\lra 2240_{10}$;

$[D_8(a_3)]\lra 1400_{11}$; $[A_7+A_1]\lra 4536_{13}$; $[E_8(a_8)]\lra 4480_{16}$;

$[E_7]\lra 84_4$; $[E_7(a_1)]\lra 567_6$; $[E_7(a_2)]\lra1344_8$; $[E_7(a_3)]\lra 2268_{10}$;

$[E_7(a_4)]\lra7168_{17}$;  $[A'_7]\lra 6075_{14}$;

$[A''_7]\lra175_{12}$;

$[A_6+A_1]\lra 2835_{14}$;

$[D_5+A_2]\lra 840_{14}$; $[D_5(a_1)+A_2]\lra1344_{19}$;

$[D_7]\lra 400_7$; $[D_7(a_1)]\lra 1050_{10}$; $[D_7(a_2)]\lra 4200_{12}$; $[D_4+A_3]\lra4200_{21}$;    

$[E_6+A_1]\lra 448_9$; $[E_6(a_1)+A_1]\lra 4096_{11}$; $[E_6(a_2)+A_1]\lra3150_{18}$;

$[A_4+A_3]\lra420_{20}$;

$[A_4+A_2+A_1]\lra2835_{22}$;

$[E_6]\lra 525_{12}$; $[E_6(a_1)]\lra 2800_{13}$; $[E_6(a_2)]\lra5600_{21}$;     

$[A_6]\lra4200_{15}$;      

$[D_6]\lra 972_{12}$; $[D_6(a_1)]\lra5600_{15}$; $[D_6(a_2)]\lra4200_{18}$; $[(2A_3)']\lra3240_{31}$;

$[D_5+A_1]\lra3200_{16}$; $[D_5(a_1)+A_1]\lra6075_{22}$;      

$[(A_5+A_1)'']\lra2016_{19}$;  

$[A_4+A_2]\lra4536_{23}$;

$[A_4+2A_1]\lra4200_{24}$;

$[D_4+A_2]\lra168_{24}$; $[D_4(a_1)+A_2]\lra2240_{28}$;

$[(2A_3)'']\lra840_{26}$;

$[A_3+A_2+A_1]\lra1400_{29}$;

$[2A_2+2A_1]\lra175_{36}$;

$[A_5]\lra 3200_{22}$     

$[D_5]\lra2100_{20}$; $[D_5(a_1)]\lra2800_{25}$;   

$[A_4+A_1]\lra4096_{26}$;

$[D_4+A_1]\lra700_{28}$; $[D_4(a_1)+A_1]\lra1400_{32}$;

$[A_3+A_2]\lra972_{32}$;

$[(A_3+2A_1)'']\lra1050_{34}$;

$[2A_2+A_1]\lra448_{39}$;

$[A_2+3A_1]\lra400_{43}$; 

$[D_4]\lra525_{36}$; $[D_4(a_1)]\lra1400_{37}$;

$[A_4]\lra2268_{30}$;

$[A_3+A_1]\lra1344_{38}$;

$[2A_2]\lra700_{42}$;

$[A_2+2A_1]\lra560_{47}$;

$[(4A_1)'']\lra 50_{56}$; 

$[A_3]\lra567_{46}$;

$[A_2+A_1]\lra210_{52}$;

$[3A_1]\lra84_{64}$;   

$[A_2]\lra 112_{63}$;

$[2A_1]\lra 35_{74}$;

$[A_1]\lra8_{91}$;     

$[A_0]\lra1_{120}$.

Here the first fourteen items are in $CL_{dist}(W)$. 
There are two conjugacy classes $A_7$; one of them, $A_7'$, comes from $W_J$ of type $E_7$.
There are two conjugacy classes $2A_3$; one of them, $2A_3'$, comes from $W_J$ of type $D_6$.
There are two conjugacy classes $A_5+A_1$; one of them, $(A_5+A_1)'$, comes from $W_J$ of
type $E_6$.
There are two conjugacy classes $A_3+2A_1$; one of them, $(A_3+2A_1)'$, comes from $W_J$ of
type $D_5$.
There are two conjugacy classes $4A_1$; one of them, $(4A_1)'$, comes from $W_J$ of
type $D_4$.

Note that $3240_9,4536_{11}$ are in $\SS^{(r)}_{dist}(W)$ for a single $r$ namely $r=2$;
$1400_{11}$ is in $\SS^{(r)}_{dist}(W)$ for a single $r$ namely $r=3$.

\head 2. A graph structure on $cl(W)$\endhead
\subhead 2.1 \endsubhead
Let $cl_{ell}(W)=\{C\in cl(W);m(C)=0\}$ (elliptic conjugacy classes). Let
$cl_{s-ell}(W)=\{C\in cl(W);m(C)=1\}$ (sub-elliptic conjugacy classes).
We shall now define a subset $\ce_{ell}(W)$ of $cl_{ell}(W)\T cl_{s-ell}(W)$. 
If $W$ is a product $W_1\T W_2$ of two Weyl group and if $\ce_{ell}(W_1),\ce_{ell}(W_2)$
are already defined, then $\ce_{ell}(W)$ consists of $(C_1\T C_2,C'_1\T C'_2)$ where either
$(C_1,C'_1)\in\ce_{ell}(W_1),C_2=C'_2\in cl_{ell}(W_2)$ or
$(C_2,C'_2)\in\ce_{ell}(W_2),C_1=C'_1\in cl_{ell}(W_1)$.
In this way we see that it is enough to define $\ce_{ell}(W)$ when $W$ is irreducible.

When $W$ is of type $A$, we set $\ce_{ell}(W)=\emp$. 

Assume that $W$ is of type $B_n,n\ge2$. With the identification 1.3(b), we define
$\ce_{ell}(W)$ to be the set of all
$((r_*,p_*),(r'_*,p'_*))\in(R\T P)^{2n}\T(R\T P)^{2n}$
such that for some $2t\in2\NN+2$ the following holds:

(a) $r'_*$ is obtained from $r_*$ by removing two consecutive terms equal to $2t$ 
and $2t$ appears at least once in $r'_*$;

(b) $p_*=(0,0,\do)$, $p'_*=(2t,2t,0,0,\do)$.

Assume that $W$ is of type $D_n,n\ge4$. With the identification 1.7(d), we define
$\ce_{ell}(W)$ to be the set of all
$((r_*,p_*),(r'_*,p'_*))\in (R_{ev}\T P)^{2n}\T(R_{ev}\T P)^{2n}$
such that for some $2t\in2\NN+2$, (a) and (b) hold.

For conjugacy classes in $W$ of exceptional type
we use the notation of Carter \cite{C72} except that
we sometimes write $D_4+A_3$ for what Carter denotes by $D_6(a_2)+A_1$.

\subhead 2.2 \endsubhead
If $W$ is of type $G_2$ we set $\ce_{ell}(W)=\emp$.

\subhead 2.3 \endsubhead
If $W$ is of type $F_4$, $\ce_{ell}(W)$ consists of:

$([D_4],[B_3])$; $([C_3+A_1],[C_3])$; $([4A_1],[3A_1])$.

\subhead 2.4 \endsubhead
If $W$ is of type $E_6$, $\ce_{ell}(W)$ consists of:

$([A_5+A_1],[A_5])$; $([3A_2],[2A_2+A_1])$.

\subhead 2.5 \endsubhead
If $W$ is of type $E_7$, $\ce_{ell}(W)$ consists of:

$([D_6+A_1],[D_6])$; $([D_6(a_2)+A_1],[D_6(a_2)])$; $([D_4+3A_1],[D_4+2A_1])$;

$([A_5+A_2],[(A_5+A_1)'])$; $([2A_3+A_1],[A_3+A_2+A_1])$; $([7A_1],[6A_1])$.

\subhead 2.6 \endsubhead
If $W$ is of type $E_8$, $\ce_{ell}(W)$ consists of:

$([E_7+A_1],[E_7])$; $([E_7(a_2)+A_1],[E_7(a_2)])$; $([E_7(a_4)+A_1],[E_7(a_4)])$;

$([E_6+A_2],[E_6+A_1])$; $([E_6(a_2)+A_2],[E_6(a_2)+A_1])$;

$([D_8(a_1)],[D_7])$; $([D_6+2A_1],[D_6+A_1])$;

$([D_5(a_1)+A_3],[D_5(a_1)+A_2])$; $([2D_4],[D_4+A_3])$;   

$([2D_4(a_1)],[D_4(a_1)+A_3])$; $([D_4+4A_1],[D_4+3A_1])$;

$([A_5+A_2+A_1],[A_5+A_2])$; $([A_5+A_2+A_1],[A_5+2A_1])$;

$([2A_4],[A_4+A_3])$;

$([2A_3+2A_2],[A_3+A_2+2A_1])$;

$([2A_3+2A_2],[2A_3+A_1])$;

$([4A_2],[3A_2+A_1])$;

$([8A_1],[7A_1])$.

\subhead 2.7 \endsubhead
We define a subset $\ce(W)$ of $cl(W)\T cl(W)$ as follows.
We say that $(C,C')\in cl(W)\T cl(W)$ is in $\ce(W)$ if there exists $J\sub I$
and $(C_1,C'_1)\in \ce_{ell}(W_J)$ such that $C=\r_J(C_1),C'=\r_J(C'_1)$, ($\r_J$ as in 0.6).
Note that $\ce_{ell}(W)\sub\ce(W)$ and that if $(C,C')\in \ce(W)$ then $m(C')=m(C)+1$.
We can regard $\ce(W)$ as the set of edges of a graph with vertices $cl(W)$.
This graph is oriented: the edge $(C,C')$ is oriented from $C$ to $C'$.
From the results of \cite{L15} one can verify that:

(a) {\it the strata of $cl(W)$ (or $W$) are exactly the connected components of this graph.}
\nl
Now each stratum of $cl(W)$ (or $W$) can be viewed as the set of vertices of an oriented graph
(restriction of the graph above to the stratum). From the results of \cite{L15} one can
verify the following strengthening of 0.4(a) and 0.4(b):

(b) {\it this oriented graph is a product of finitely many Coxeter graphs
of type $A$ (with the usual orientation).}
\nl
Recall that an oriented Coxeter graph of type $A$ is of the form
$$\bul\to\bul\to\do\to\bul.$$

For example, if $W$ is of type $E_8$, then the graph attached to a stratum of $cl(W)$ (or $W$)
is of one of the types $A_5,A_4,A_3,A_2\T A_2,A_2,A_1$. (Type $A_5$ appears for a unique
stratum; type $A_2\T A_2$ appears for two strata.)
If $W$ is of type $E_7$, then the graph attached to a stratum of $cl(W)$ (or $W$)
is of one of the types $A_4,A_3,A_2,A_1$. (Type $A_4$ appears for a unique
stratum.) If $W$ is of type $E_6$, then the graph attached to a stratum of $cl(W)$ (or $W$)
is of one of the types $A_2,A_1$.
If $W$ is of type $F_4$, then the graph attached to a stratum of $cl(W)$ (or $W$)
is of one of the types $A_4,A_2,A_1$. (Type $A_4$ appears for a unique stratum.)
If $W$ is of type $G_2$. then the graph attached to a stratum of $cl(W)$ (or $W$)
is of type $A_1$. 

\head 3. Complements\endhead
\subhead 3.1 \endsubhead
Let $r\in\{0\}\cup\cp$. We state two properties similar to 0.4(a),(b).

(a) For any $E\in\SS^{(r)}(W)$  there is a unique $C^{(r)}_E\in(\Ph^{(r)})\i(E)$ which is as elliptic as
possible, that is $m(C_E^{(r)})\le m(C)$ for any $C\in(\Ph^{(r)})\i(E)$;

(b) For any $E\in\SS^{(r)}(W)$  there is a unique $C'{}^{(r)}_E\in(\Ph^{(r)})\i(E)$ which is as non-elliptic as
possible, that is $m(C'{}^{(r)}_E)\le m(C)$ for any $C\in (\Ph^{(r)})\i(E)$.
\nl
Now (a) is proved in \cite{L12}; the proof of (b) is similar.
Let $CL^{(r)}(W)$ be the image of the map $E\m C'_E$, $\SS^{(r)}(W)@>>>cl(W)$. Note that
$\Ph^{(r)}$ restricts to a bijection $CL^{(r)}(W)@>\si>>\SS^{(r)}(W)$.
Let $CL_{dist}^{(r)}(W)$ be the subset of $CL^{(r)}(W)$ corresponding to $\SS^{(r)}_{disc}\sub\SS^{(r)}(W)$
under this bijection. We have the following analogue of 0.4(c).

(c) {\it Let $C\in CL^{(r)}(W)$. We have $C\in CL_{dist}^{(r)}(W)$ if and only if $C$ is elliptic.}
\nl
Note that $CL_{dist}(W)=\cup_rCL_{dist}^{(r)}(W)$.

\subhead 3.2 \endsubhead
For any semisimple element $s\in G=G^{(0)}$ let $E_s=j_{W_s}^W(\text{sign})$ ($j$-induction)
where $W_s$ is the Weyl group of the connected centralizer of $s$ viewed as a subgroup of $W$.
We have $E_s\in\SS(W)$; the subset of $\SS(W)$ formed by the $E_s$ for various $s$ as above is
denoted by $\Irr_{ss}(W)$.

Now let $E\in\SS_{dist}(W)$ and let $C\in CL_{dist}(W)$ be the corresponding conjugacy class.
According to 0.5(a) we have $E=\Ph^{(0)}(C)$
hence $E=\io^{(0)}(\g)$ where $\g=\ps^{(0)}(C)\in\cu^{(0)}$.
(The elements of $\g$ need not be distinguished.)
From \cite{L11b} there exists a semisimple element $s\in G^{(0)}$ such that
$E=E_s$; thus, we have $\SS_{dist}(W)\sub\Irr_{ss}(W)$.
(We can assume that $G$ is almost simple; then
the statement in the previous sentence holds for $E=\Ph^{(0)}(C)$ for any
elliptic $C\in cl(W)$ with a single exception
in type $E_8$; but that exception is not distinguished, see 1.15.)
Note that $s$ belongs to the stratum of $G^{(0)}$ that contains $\g$.

Thus, to $E\in\SS_{dist}(W)$ (or $C\in CL_{dist}(W)$)
one can associate a collection of reflection subgroups $W_s$
(for various $s$ as above).

Now, \cite{L11a, 4.4(b)} implies that the minimum length of
an element in $C\in CL_{dist}(W)$ is equal to the dimension of the centralizer of $s$ (as above)
in $G$ modulo the centre of $G$.

\subhead 3.3\endsubhead
In the 3.4-3.10 we describe explicitly a correspondence

(a) $C\m W_s$
\nl
which to any $C\in CL_{dist}(W)$ associates a reflection subgroup $W_s$ of $W$ (up to
conjugacy) as in 3.2 (assuming that $W$ is irreducible of type $\ne A$); if $W$ is of type
$A$, then $C$ is the Coxeter class and the corresponding $W_s$ is $\{1\}$.

\subhead 3.4\endsubhead
We now assume that $G=Sp_{2n}(\CC),n\ge2$ or that $G=SO_{2n+1}(\CC),n\ge2$.
According to 1.4, 1.6, $CL_{dist}(W)$
can be identified with the set of pairs $(r_*,p_*)\in(R\T P)^{2n}$ such that $p_*=(0,0,\do)$
and $r_*$ is a sequence $r_1\ge r_2\ge r_3\ge\do\ge r_\s$ of even integers $>0$ without
two consecutive equalities. For such $(r_*,p_*)$ we define
$k=r_1/2$. We define $\br_1\ge\br_2\ge\do\ge\br_k>0$
by $\br_i=\sha(t\in[1,\s];r_t/2\ge i)$ for $i\in[1,k]$; note that $\br_1=\s$ and
$\br_1+\br_2+\do+\br_k=n$.
Using 0.5(a) and \cite{L11b},
we see that the reflection subgroup $W_s$ corresponding to $(r_*,p_*)$
is a product of symmetric groups
$$S_{\br_k}\T S_{\br_{k-1}}\T\do\T S_{\br_1}$$
(a subgroup of $S_n$ which is itself naturally a subgroup of $W$ of the form $W_J$ for some
$J$).

We see that $CL_{dist}(W)$ is in natural bijection with the set of sequences
$\br_1\ge\br_2\ge\br_3\ge\do$ of integers $\ge0$ such that
$\br_1+\br_2+\do=n$ and $\br_i-\br_{i+1}\le2$ for all $i\ge1$.

\subhead 3.5\endsubhead
We now assume that $G=S0_{2n}(\CC),n\ge4$. According to 1.8, $CL_{dist}(W)$
can be identified with the set of pairs $(r_*,p_*)\in(R\T P)^{2n}$ such that $p_*=(0,0,\do)$
and $r_*$ is a sequence $r_1\ge r_2\ge r_3\ge\do\ge r_\s$ of even integers $>0$ without
two consecutive equalities with $\s$ even. For such $(r_*,p_*)$ we define
$\br_1\ge\br_2\ge\do\ge\br_k>0$ as in 3.4. Note that $\br_1$ is even and $\ge2$.
Using 0.5(a) and \cite{L11b},
we see that the reflection subgroup $W_s$ corresponding to $(r_*,p_*)$ is
of the form
$$S_{\br_k}\T S_{\br_{k-1}}\T\do\T S_{\br_2}\T W_{D_{\br_1/2}}\T W_{D_{\br_1/2}}.$$
Here $W_{D_m}$ denotes a Weyl group of type $D_m$ (for $m=1$ this is taken to be $\{1\}$;
for $m=2$ this taken to be $S_2\T S_2$). We view
$$S_{\br_k}\T S_{\br_{k-1}}\T\do\T S_{\br_2}\T W_{D_{\br_1}}$$
as a subgroup of $W$ of the form $W_J$ for some $J$ and $W_{D_m}\T W_{D_m}$ as a
 reflection subgroup of $W_{D_{2m}}$ in the standard way.

We see that $CL_{dist}(W)$ is in natural bijection with the set of sequences
$\br_1\ge\br_2\ge\br_3\ge\do$ of integers $\ge0$ such that
$\br_1+\br_2+\do=n$, $\br_1$ is even and $\br_i-\br_{i+1}\le2$ for all $i\ge1$.

\subhead 3.6\endsubhead
Until the end of 3.10, $W$ is of exceptional type.
In each case the reflection subgroup $W_s$ attached by 3.3(a)
to $C\in CL_{dist}(W)$ is specified by its type. (We use 0.5(a) and \cite{L11b}.)

If $W$ is of type $G_2$, the correspondence 3.3(a) is:

$[G_2]\m A_0$; $[A_2]\m A_1$; $[A_1+\tA_1]\m 2A_1$.

Note that the group $W_s$ is of the form $W_J$ for some $J\sub I$
except for the last case.

\subhead 3.7\endsubhead
If $W$ is of type $F_4$, the correspondence 3.3(a) is:

$[F_4]\m A_0$; $[B_4]\m A_1$; $[F_4(a_1)]\m 2A_1$; $[D_4(a_1)]\m A_2+A_1$;

$[A_3+\tA_1]\m B_2+A_1$; $[\tA_2+\tA_2]\m A_2+A_2$.

Note that the group $W_s$ is of the form $W_J$ for some $J\sub I$
except for the last two cases.

\subhead 3.8\endsubhead
If $W$ is of type $E_6$, the correspondence 3.3(a) is:

$[E_6]\m A_0$; $[E_6(a_1)]\m A_1$; $[E_6(a_2)]\m 3A_1$.

Note that the group $W_s$ is of the form $W_J$ for some $J\sub I$.

\subhead 3.9\endsubhead
If $W$ is of type $E_7$, the correspondence 3.3(a) is:

$[E_7]\m A_0$; $[E_7(a_1)]\m A_1$; $[E_7(a_2)]\m 2A_1$; $[E_7(a_3)]\m 3A_1$;

$[A_7]\m A_2+2A_1$; $[E_7(a_4)]\m 2A_2+A_1$.

Note that the group $W_s$ is of the form $W_J$ for some $J\sub I$.

\subhead 3.10\endsubhead
If $W$ is of type $E_8$, the correspondence 3.3(a) is:

$[E_8]\m A_0$; $[E_8(a_1)]\m A_1$; $[E_8(a_2)]\m 2A_1$; $[E_8(a_4)]\m 3A_1$;

$[E_8(a_5)]\m 4A_1$; $[D_8]\m A_2+2A_1$; $[E_8(a_3)]\m A_2+3A_1$;

$[E_8(a_7)]\m 2A_2+A_1$; $[E_8(a_6)]\m 2A_2+2A_1$;

$[A_8]\m A_3+A_2+A_1$; $[E_8(a_8)]\m A_4+A_3$;

$[D_8(a_2)]\m A_3+3A_1$; $[D_8(a_3)]\m A_3+A_2+2A_1$; $[A_7+A_1]\m 2A_3+A_1$.

Note that the group $W_s$ is of the form $W_J$ for some $J\sub I$
except for the last three cases.

\subhead 3.11\endsubhead
Let $K$ be a maximal compact subgroup of $G^{(0)}$. The following result was stated in \cite{L21, 5.2}:

(a) {\it Let $X$ be a stratum  of $G^{(0)}$ and let $E$ be
the corresponding element of $\SS(W)$. We have $X\cap K\ne\emp$ if and only if $E\in\Irr_{ss}(W)$.}
\nl
By the results in 3.2, we have $E\in \Irr_{ss}(W)$ if and only if $X$
contains a semisimple element of $G^{(0)}$. This last condition is clearly satisfied when
$X\cap K\ne\emp$. Conversely, assume that $X$ contains a semisimple element $s$ of $G^{(0)}$. It is
well known that we can find $s'\in K$ such that the connected centralizers of $s$ and $s'$ are
conjugate. It follows that $E_s=E_{s'}$ hence $s,s'$ belongs to the same stratum. Since $s\in X$ we
have $s'\in X$ so that $X\cap K\ne\emp$. This proves (a).

We show:

(b) {\it Let $X,E$ be as in (a). If $E$ is distinguished then $X\cap K\ne\emp$.}
\nl
By (a) it is enough to show that $E\in\Irr_{ss}(W)$.
This follows from $\SS_{dist}(W)\sub\Irr_{ss}(W)$, see 3.2.

{\it Errata to \cite{L15}}.

p.355, line  containing $378_{14}$: replace $[A_3+A_2]$ by $[2A_3]$.

p.356, line containing $3200_{22}$: replace $[(A_5+A_1)']$ by $[(A_5+A_1)',A_5]$.

\enddocument